\documentclass[12pt]{article}
\usepackage{amsmath,amsthm,amsfonts,latexsym,amssymb}
\def\idp{\mathfrak{p}}  
\def\idP{\mathfrak{P}} 
\def\ida{\mathfrak{a}}  
\def\idb{\mathfrak{b}}  
\def\idc{\mathfrak{c}}  
\def\Q{{\mathbb Q}}     
\def\Z{{\mathbb Z}}     
\def\E{{\mathcal{E}}}
\def\O{{\mathcal{O}}}
\def\S{{\mathcal{S}}} 
\def\al{{\alpha}} 
\def\be{{\beta}}  
\def\ga{{\gamma}}  
\def\de{{\delta}}  
\def\th{{\theta}}  
\def\ep{{\epsilon}}  
\def\om{{\omega}}  

\begin{document}
\title{On the equation $Y^2 = X^6 + k$}

\author{A.~Bremner\thanks{Department of Mathematics, Arizona State University, 
Tempe AZ, USA, e-mail: bremner@asu.edu, 
http://math.asu.edu/\~{}andrew/} \and
N.~Tzanakis\thanks{Department of Mathematics, University of Crete, Iraklion, 
Greece, e-mail: tzanakis@math.uoc.gr, http://www.math.uoc.gr/\~{}tzanakis}
}
\date{January 20, 2010}
\maketitle

\begin{center}
Dedicated with respect and admiration to Professor Paulo Ribenboim \\ 
on the occasion of his eightieth birthday
\end{center}
%
\begin{abstract}
We find explicitly all rational solutions of the title equation for all
integers $k$ in the range $|k|\leq 50$ except for $k=-47,-39$. For the
solution, a variety of methods is applied, which, depending on $k$, may range
from elementary, such as divisibility and congruence considerations, to
elliptic Chabauty techniques and highly technical computations in algebraic
number fields, or a combination thereof. For certain sets of values of $k$ we
can propose a more or less uniform method of solution, which might be applied
successfully for quite a number of cases of $k$, even beyond the above range.
It turns out, however, that in the range considered, six really challenging
cases have to be dealt with individually, namely $k = 15,43,-11,-15,-39,-47$.
More than half of the paper is devoted to the solution of the title equation
for the first four of these values. For the last two values the solution of 
the equation, at present, has resisted all our efforts. The case with these
six values of $k$ shows that one cannot expect a general method of solution
which could be applied, even in principle, for {\em every} value of $k$.
A summary of our results is shown at the end of the paper.
\end{abstract}

%
\section{Introduction} \label{introduction}
%
For a fixed integer $k$, there is a vast literature devoted to the integer solutions 
of the Diophantine equation $Y^2=X^3+k$ (Mordell's equation). It is well known 
that the equation has only finitely many integer solutions, and many individuals 
over many years investigated this problem, starting with solving
the equation for particular values of $k$. There now exists a very satisfactory 
uniform method (the elliptic logarithm method, developed independently by
Stroeker and Tzanakis~\cite{ST} and Gebel, Peth\H{o}, Zimmer~\cite{GPZ1}) for explicitly 
computing all integer solutions on elliptic curves.  See Gebel, Peth\H{o}, Zimmer~\cite{GPZ2}
for a very successful application of this method to computing all integer points 
on Mordell's equation for a given $k$. 
When we consider {\it rational} solutions to Mordell's
equation, then the whole theory of elliptic curves comes into play, and many
questions remain unresolved, which is not however the focus here. In this paper 
we shall study rational solutions of the naturally arising Diophantine equation 
\begin{equation}
\label{kcurve}
C_k:Y^2=X^6+k
\end{equation}
for a fixed integer $k \neq 0$, free of sixth powers. There are two rational 
points at infinity on $C_k$ corresponding to $Y=\pm X^3$, and henceforth we shall 
consider only finite rational
points. Since $C_k$ represents a curve of genus 2, Faltings's theorem (``Mordell's 
Conjecture")~\cite{F}, shows that $C_k(\Q)$ is finite. 
It is well known that both Falting's proof and the subsequent one by 
P.~Vojta~\cite{V} are non-effective (an effectively computable bound for the number of
rational points is accomplished), hence they cannot provide us with a
practical method for the explicit determination of all rational points on $C_k$.
Thus determination of rational points on $C_k$ is a good challenge. We attack this 
challenge by a series of approaches that start with elementary ideas and progress to
being more technical; we illustrate the ideas by applying them to $k$ in the range
$|k| \leq 50$. \\ \\
If the rank of the Jacobian $J_k$ of $C_k$ is at most equal to 1, then there are effective 
``Chabauty" techniques that allow (at least in principle, but not provably) determination 
of $C_k(\Q)$. However, the Jacobian $J_k$ is isogenous to the product of the two 
elliptic curves
\begin{equation}
\label{E1E2}
\E_1: Y_1^2=X_1^3 + k, \qquad \E_2: Y_2^2=X_2^3+k^2,
\end{equation}
with maps from $C_k$ to $\E_i$, $i=1,2$, given by
\[ (X_1,Y_1)=(X^2,Y), \qquad  (X_2,Y_2)=(\frac{k}{X^2}, \frac{k Y}{X^3}). \]
Thus the rank of $J_k$ is the sum of the ranks of $\E_1$ and $\E_2$. If $J_k$ has rank 
at most 1, then there is no need to use Chabauty techniques, for in such a case 
the rank of $\E_1$ or $\E_2$ is 0, and either
there are no finite rational points on (\ref{kcurve}), or there are only
the obvious ones arising from non-zero torsion points on $\E_i$ (which occur 
for $k$ a perfect square or perfect cube).\\ \\
The interesting cases therefore arise when both $\E_1$ and $\E_2$ have positive ranks.
We know of only one example in the literature where the set of rational points is 
determined upon such a curve: Coombes and Grant~\cite{CG} in an example show that
the only rational points on the curve $y^2=x^6-972$ (with covering elliptic curves
both of rank 1) are the points at infinity. This will also follow from the 
elementary arguments of section (\ref{elementary args}) below. \\
To deal with (\ref{kcurve}), we shall construct maps to $C_k$ from curves having 
ternary equations of type 
$a x^6+b y^6=c z^3$. This is the particular form of a very natural generalized Fermat curve.
Several important papers recently have been devoted to generalized Fermat equations,
and papers of Bennett \& Skinner~\cite{BS} and Darmon \& Granville~\cite{DG}
seem applicable when at least one of the exponents in the equation is large. However,
the specific equations $a x^6+b y^6=c z^3$ (and $a x^6+b y^6=c z^2$) seem
poorly treated, indeed, essentially absent, in the literature. We can find one
paper by Dem'yanenko~\cite{Dem}, who inter alia discusses the equations
$x^6+y^6=a z^2$ and $x^6+y^6=a z^3$. 
This dearth is curious. For several decades
in the early twentieth century, number theorists were attracted to showing the
impossibility of specific Diophantine equations by elementary methods; we have
here a class of such equations where impossibility may often be shown by an
argument involving an elliptic curve of rank 0, in other words, by a classical
elementary infinite descent argument. This class of equations seems to have been 
overlooked
by earlier researchers. 
In any event, for a certain class of $k$, these curves $a x^6+b y^6=c z^3$ of genus 4 
may be shown by elementary means to have no non-trivial rational points. This
furnishes explicitly all solutions to (\ref{kcurve}) 
for the corresponding $k$; see section (\ref{generalities}). Within this class, there 
are values of $k$ for which the descent argument, while elementary, is not straightforward, 
and hence deserves
some discussion; see section (\ref{elementary}). Some values of $k$ may also
be eliminated using the alternative elementary ideas of section (\ref{congruences}). \\ \\
For many values of $k$, the above descent is inconclusive, and we turn to working
over $K=\Q(\theta)$, with $\theta^3=k$. There arises the problem of solving a
number of Diophantine problems of the following type: for a certain elliptic curve
$E$ defined over $K$ such that the rank of $E(K)$ is at most 2, explicitly determine
those points $(x,y) \in E(K)$ satisfying a certain ``rationality condition"
$q(x,y) \in \Q$ where $q(X,Y) \in K(X,Y)$. Problems of this type are amenable
to the elliptic Chabauty method as implemented in a number of routines
in Magma~\cite{Mag}. But we should stress that the application of Magma
is in general not a matter of button pushing, and is sometimes far from automatic;
we note later some of the problems that can arise.  For
a certain class of $k$, the elliptic Chabauty method applies more or less directly;
and we consider the corresponding equation (\ref{kcurve}) to be interesting
but whose solution is ``relatively standard''. So few details are given;
see section (\ref{Chabauty}). \\ \\
For several values of $k$ in our considered range, the elliptic Chabauty method
as applied above does not succeed, as we explain in the first
paragraph of section (\ref{difficultcases}).
This happens for $k=-47$, $-39$, $-15$, $-11$, $15$, $43$.
We surmounted obstacles that arise by means of an alternative approach for the four cases
$k=15$, $43$, $-11$, $-15$ (given in decreasing order of difficulty, the first two being roughly 
comparable). We give details on
solving (\ref{kcurve}) for $k=15$ in section (\ref{k=15}), with brief details for
$k=43$ in section (\ref{k=43}); and for $k=-11,-15$ in sections (\ref{k=-11}), (\ref{k=-15}).
These are the most technical sections of the paper.\\
We have been unable to develop a uniform method for the $k$ of this section,
though believe the ideas will be useful for solving (\ref{kcurve})
for other specific ``difficult" values of $k$.\\ \\
The remaining values $k=-47,-39$ in our range have resisted all our attacks,
These values give no indication of the difficulties that arise when trying to solve 
(\ref{kcurve}), difficulties that may not be foreseen until after many hours of work. 
Finding all rational points on (\ref{kcurve}) for these intractable values of $k$ 
is indeed a challenging Diophantine problem.
%
\section{General discussion} \label{generalities}
%
In this section we give some general discussion, the emphasis being on elementary arguments which often 
suffice to provide the complete solution of equation (\ref{kcurve}). Our arguments will then be applied
to solve a good number of equations (\ref{kcurve}) in the range $|k|\leq 50$. 
Of course, as one should expect, the equation (\ref{kcurve}) is not susceptible
to elementary treatment for all $k$, not even within the small range we consider. 
More technical approaches are given in sections (\ref{Chabauty}), (\ref{difficultcases}), and (\ref{defeat}) of this paper.
%
\subsection{When at least one $\E_i$ has zero rank} \label{zero rank cE}
If either of the curves $\E_1$ or $\E_2$ at (\ref{E1E2}) has rank zero, then
the determination of rational points on $C_k$ is trivial.  For if $\E_1$ has rank 0, 
the only possible rational points on $\E_1$ are torsion points, and using, for example, 
Cassels~\cite{Ca}, Theorem p.52, we have on $\E_1$ for $k=-432$ the points $(12,\pm 36)$; 
for $k=1$ the points $(2,\pm 3)$, $(0,\pm 1)$, $(-1,0)$; for $k=D^3$, $D \neq 1$, the point
$(-D,0)$; and for $k=B^2$, $B \neq 1$, the points $(0,\pm B)$. In other cases, the 
torsion group is trivial. It follows that the only (finite) points on $C_k$ occur when 
$k=-1$, with points $(\pm 1, 0)$; and when $k=B^2$, with points $(0,\pm B)$. \\
Similarly, if $\E_2$ has rank 0, the same Theorem shows that no (finite)
points arise on $C_k$. \\
In the range $|k|\leq 50$ at least one curve $\E_i$ has rank 0 when
\begin{align*}
  k= & \pm 1,\pm 2, -3,\pm 4, \pm 5, \pm 6, \pm 7,\pm 8, \pm 9, -10, \pm 12, 13,
       \pm 14, \pm 16, -17, \pm 18, \pm 19,\pm 20, 21, \pm 22, \\
     &  \pm 23, -24, 25, \pm 26, \pm 27, 29, \pm 30, 
     -31, \pm 32, \pm 33, \pm 34, -36, \pm 37, \pm 38, \pm 40, \pm 41, \pm 42, \pm 44, \\
    & \; 45, -46, 49, \pm 50\:.
\end{align*}
So, for the above $k$'s equation (\ref{kcurve}) is immediately solved. 
%
\subsection{Elementary arguments} \label{elementary args}
Henceforth we shall assume that the ranks of both $\E_1$ and $\E_2$ are positive. 
In the range $|k|\leq 50$, this is the case for for the following values:
\begin{align}
 k = & -49, -48, -47,-45,-43, -39, -35, -29, -28, -25, -21, -15, -13, -11,  \nonumber \\
   & 3, 10, 11, 15, 17, 24, 28, 31, 35, 36, 39, 43, 46, 47, 48\:.  
                                        \label{k-sieve 1}
\end{align}
Now, finding all finite points on (\ref{kcurve}) is equivalent to solving
\begin{equation}
 x^6+k y^6 = z^2, \quad x,y,z \in \Z, \qquad y\neq 0,\;(x,y)=1, \label{xyeq}
\end{equation}
where $X=x/y$ and $Y=z/y^3$. 
The factorization $(z+x^3)(z-x^3)=k y^6$ of (\ref{xyeq}) is easily seen to lead 
to a number of equations of type 
\begin{equation} \label{ABsextic}
 A y_1^6+B y_2^6=y_3^3, \;\; A,B \in \Z\;; 
\end{equation}
and the curve (\ref{ABsextic}) has maps to elliptic curves as follows, where the 
notation $E:=[a,b,c,d,e]$ means that $E$ is the elliptic curve with Weierstrass coefficients 
$a,b,c,d,e$:
\begin{equation} \label{e1e2}
(\frac{A y_3}{y_2^2}, \frac{A^2 y_1^3}{y_2^3})\;\mbox{ on }\;E_1:=[0,0,0,0,-A^3 B], \quad 
(\frac{B y_3}{y_1^2}, \frac{B^2 y_2^3}{y_1^3})\;\mbox{ on }\;E_2:=[0,0,0,0,-A B^3]; 
\end{equation}
and
\begin{equation} \label{e3}
( \frac{y_3^6-A B y_1^6 y_2^6}{y_1^4 y_2^4 y_3^2},\frac{(A y_1^6-B y_2^6)(2A y_1^6+B y_2^6)(A y_1^6+2B y_2^6)}{2 y_1^6 y_2^6 y_3^3} )\;\mbox{ on }\;E_3:=[0,0,0,0,-\frac{27}{4} A^2 B^2].
\end{equation}
If the rank of $E_1$ or $E_2$ is 0, then the solution of (\ref{ABsextic}) is
immediate. (The curve $E_3$ is in fact 3-isogenous to the curve $\E_2$, so by previous assumption
has positive rank). \\

According as to the possibilities for $k\bmod{\,4}$, four distinct types of 
equation (\ref{ABsextic}) arise, which we denote by Types I, II, III, IV, as follows.

\noindent
\textbullet\; If $k\equiv1\pmod{2}$.
\[ \mathrm{I}:\quad 2x^3=\de_1y_1^6+\de_2y_2^6   \quad y_1y_2\neq 0\,, \qquad (A,B)=(4\de_1,4\de_2)\,, \]
\[y=y_1y_2,\quad  \de_1\de_2=-k,\; \de_2>0,\quad (\frac{\de_1}{\de}y_1,\frac{\de_2}{\de}y_2)=1,\,
\: \de=(\de_1,\de_2),\quad   (2x,y_1y_2)=1 \:.\]
%
\vspace*{0.2in}
\[ \mathrm{II}:\quad x^3=16\de_1y_1^6+\de_2y_2^6   \quad y_1y_2\neq 0 \,, \qquad (A,B)=(16\de_1,\de_2)\,, \]
\[y=2y_1y_2,\quad  \de_1\de_2=-k,\; \de_2>0,\quad (\frac{2\de_1}{\de}y_1,\frac{\de_2}{\de}y_2)=1,\,
\: \de=(\de_1,\de_2),\quad   (x,2y_1y_2)=1 \:.\]
%
\textbullet\; If $k\equiv0\pmod{4}$.
\[ \mathrm{III}:\quad x^3=\de_1y_1^6+\de_2y_2^6   \quad y_1y_2\neq 0 \,, \qquad (A,B)=(\de_1,\de_2)\,, \]
\[y=y_1y_2,\quad  \de_1\de_2=-\frac{k}{4},\; \de_2>0,\quad (\frac{\de_1}{\de}y_1,\frac{\de_2}{\de}y_2)=1,\,
\: \de=(\de_1,\de_2),\quad   (x,y_1y_2)=1 \:.\]
%
\textbullet\; If $k\equiv 2\pmod{4}$.
\[ \mathrm{IV}:\quad x^3=32\de_1y_1^6+\de_2y_2^6   \quad y_1y_2\neq 0 \,, \qquad (A,B)=(32\de_1, \de_2)\,, \]
\[y=2y_1y_2,\quad  \de_1\de_2=-\frac{k}{2},\; \de_2>0,\quad (\frac{2\de_1}{\de}y_1,\frac{\de_2}{\de}y_2)=1,\,
\: \de=(\de_1,\de_2),\quad   (x,2y_1y_2)=1 \:.\]
%
Equation (\ref{xyeq}) is thus reduced to a set $\S(k)$ of equations of type 
I, II, III, IV.
For the following values of $k$ at (\ref{k-sieve 1}) all equations in $\S(k)$ are 
either impossible or have only the obvious solutions:
\begin{equation}  \label{k-sieve 2}
 k=-49,-48,-45,-13,11,28,36,39,46,47 \:.
\end{equation}
This is proved by showing that each equation in $\S(k)$ either is impossible $\bmod{\;m}$, 
where $m\in\{7,8,9,13\}$, or corresponds to a curve $E_i$ at (\ref{e1e2}) with zero rank 
for either $i=1$ or $i=2$.

For the following values of $k$ all equations in $\S(k)$ are also impossible; however, 
for at least one element of $\S(k)$, more refined but still elementary arguments are 
required:
\begin{equation}  \label{k-sieve 3}
 k=-43,-35, -29, -25, -21, 31 \:.
\end{equation}
We choose to give in the following subsection some details of the proposed elementary 
method when $k=-35$ which is 
a most characteristic case. Here, $\S(-35)$ consists of two equations of type 
I
and four equations of type 
II.
Five out of the six equations are impossible 
as congruences modulo 7 or 8 and only the equation (of type 
II)
\begin{equation} \label{xyeq -35}
16\cdot 35y_1^6+y_2^6=x^3 \,,\quad y_1y_2\neq 0
\end{equation}
remains. This belongs to the more general class of equations
\begin{equation} \label{Deq}
 DY_1^6 +Y_2^6 =X^3\,,\quad Y_1\neq 0,\: (DY_1,Y_2)=1 \:,
\end{equation}
where $D$ is a sixth power free non-zero integer. 

It is worth noting that for a number of $D$'s, equations (\ref{Deq}) can be treated by quite
elementary means.
%
\subsubsection{An elementary approach to (\ref{Deq}) with application
to $k=-43,-35,-21,31$} \label{elementary}
%
We distinguish two cases, depending on the divisibility of $X-Y_2^2$ by 3.

Case (i): $X-Y_2^2\not\equiv 0\pmod{3}$. In this case $(X-Y_2^2,X^2+XY_2^2+Y_2^4)=1$, 
$DY_1\not\equiv 0\pmod{3}$ and
\[ X-Y_2^2=d_1Y_3^6\,,\quad X^2+XY_2^2+Y_2^4 = d_2Y_4^6 \:,\]
where
\[d_1d_2=D,\,d_2>0,\,(2d_1Y_3,d_2Y_4)=1,\,Y_1=Y_3Y_4\not\equiv 0\!\!\!\pmod{3},\, 
          (X,Y_2Y_3Y_4)=1,\,(Y_2,d_2)=1 \:.\]
Substitution of $X$ from the first equation into the second gives
\begin{equation} \label{quartic I}
 3Y_2^4+3d_1Y_2^2Y_3^6 +d_1^2Y_3^{12} =d_2Y_4^6 \:.
\end{equation}
Observe first that the following conditions are necessary for the solvability of
equation (\ref{quartic I}):
\begin{enumerate}
\item $(d_2,d_1)\equiv (1,1),(1,7),(1,8),(4,2),(4,4),(4,7),(7,1),(7,4),(7,5) \pmod{9}$ \\
\item $\left(\frac{3d_2}{p}\right) = 1 \quad\mbox{for every odd prime divisor $p$ of $d_1$}$ \\
\item $p\equiv 1\pmod{3} \quad\mbox{for every odd prime divisor $p$ of $d_2$} \:.$
\end{enumerate}
Note that in our example with $k=-35$, only equation (\ref{xyeq -35}) is left to treat,
and accordingly we take $D=16\cdot 35$.  It is straightforward to check that the first 
condition above is satisfied for no pair $(d_2,d_1)$
which means that (\ref{xyeq -35}) is impossible if $x\not\equiv y_2^2\!\!\pmod{3}$.

Case (ii): $X-Y_2^2\equiv 0\pmod{3}$. In this case let
\[ 3^{\tau}|| D \quad\mbox{and}\quad 
 \nu=\begin{cases} 5 & \mbox{if $\tau=0$} \\
                   6 &  \mbox{if $\tau=1$} \\
                   \tau-1 & \mbox{if $\tau\geq 2$}
     \end{cases}
\:.
\]
It is easy to see that
\[ X-Y_2^2=3^{\nu}d_1Y_3^6, \;\; X^2+XY_2^2+Y_2^4 = 3d_2Y_4^6, \mbox{ with } d_1d_2=3^{-\tau}D,\;\; d_2>0, \;\; Y_1=3^\frac{\nu+1-\tau}{6}Y_3Y_4,\]
where
\[ (2d_1Y_3,d_2Y_4)=1, \quad (X,3^{\nu+1-\tau}Y_2Y_3Y_4)=1, \quad (Y_2,d_2)=1 \:.\]

Substitution of $X$ from the first equation into the second gives
\begin{equation} \label{quartic II}
 Y_2^4+3^{\nu}d_1Y_2^2Y_3^6 +3^{2\nu-1}d_1^2Y_3^{12} =d_2Y_4^6 \:.
\end{equation}
Necessary conditions for the solvability of (\ref{quartic II}) are the following:
\begin{enumerate}
\item $\left(\frac{d_2}{p}\right)=1 \quad\mbox{for every odd prime divisor $p$ of $d_1$}$ \\
\item $p\equiv 1\pmod{3} \quad\mbox{for every odd prime divisor $p$ of $d_2$}\:.$ 
\end{enumerate}
In the case $k=-35$, $D=16\cdot 35$, we have $\tau=0$, $\nu=5$, and the conditions above
are satisfied by no pair $(d_2,d_1)$ except for $(d_2,d_1)=(1,16\cdot 35)$, for which
equation (\ref{quartic II}) has an obvious solution and hence cannot be excluded by 
congruence considerations. To proceed further, however, we factor (\ref{quartic II}) 
over $\Q(\om)$, where $\om^2+\om+1=0$, obtaining (for all values of $k$)
\begin{equation} \label{quartic II omega}
 Y_2^2+3^{\nu-1}d_1(2+\om)Y_3^6 =(m+n\om)(a+b\om)^6, \qquad m,n,a,b \in \Z,
\end{equation}
where
\[ m^2-mn+n^2=d_2\,,\quad a^2-ab+b^2 = |Y_4|\,,\;(a,b)=1 \:.\]
Moreover, if $\tau\neq 2$ we can assume without loss of generality that
\[ m\not\equiv 0,\; n\equiv 0 \!\! \pmod{3} \quad\mbox{and $ab$ is odd with 
$a+b\not\equiv 0 \!\!\pmod{3}$.}
\]
From (\ref{quartic II omega}), on equating coefficients of $\om$, $1$, we obtain
\begin{eqnarray}
 F_1(a,b) & = & 3^{\nu-1}d_1Y_3^6  \label{eqF1} \\
 F_2(a,b) & = & Y_2^2 \label{eqF2}
\end{eqnarray}
where
\begin{align*}
 F_1(a,b) = & na^6+6(m-n)a^5b-15ma^4b^2+20na^3b^3+15(m-n)a^2b^4-6mab^5+nb^6  
                                          \\[1mm]
F_2(a,b) = & (m-2n)a^6-6(2m-n)a^5b+15(m+n)a^4b^2+20(m-2n)a^3b^3  \\ 
          & -15(2m-n)a^2b^4+6(m+n)ab^5+(m-2n)b^6 \:.
\end{align*}

In our case $k=-35$, we have $d_2=1$, so $n=0$, and equation (\ref{eqF1}) becomes 
\[ \frac{ab(a-b)}{2}\cdot\frac{(a+b)(2a-b)(a-2b)}{2}=\pm 2^2 \cdot 3^3 \cdot 35 Y_3^6 \]
where the two factors on the left-hand side are relatively prime. Putting 
$a/b=u\in\Q$ we obtain 
\begin{equation} \label{c1c2}
 u(u-1)=2c_1w_1^3\,,\quad (u+1)(2u-1)(u-2)=2c_2w_2^3\,,\quad w_1,w_2\in\Q \:,
\end{equation}
where $c_1,c_2$ are relatively prime positive integers with $c_1c_2=140$. It is easily checked
that for every such pair $(c_1,c_2)$ at least one of the elliptic curves at (\ref{c1c2})
is of zero rank, from which we easily conclude that equation (\ref{xyeq -35}) is impossible
in case (ii).

Mutatis mutandis, for the remaining values of $k$ of this subsection, namely 
$k=-43$, $-21$, $31$ (as well as for $k=-972$, the example considered 
by Coombes and Grant~\cite{CG}), the set $\S(k)$ contains equations that are either impossible
or possess only trivial solutions. 
%
\subsubsection{Congruences on elliptic curves with application 
to $k=-29,-25$.} \label{congruences}
%
If an equation in the set $\S(k)$ is everywhere locally solvable, yet is suspected of 
having no rational solution, we can in some instances still prove impossibility
by elementary means using the following trick involving congruences 
with points on elliptic curves.  For the remaining values of $k$ in the range $|k|\leq 50$,
this applies to $k=-29,-25$ (see also $k=-15$ in section (\ref{k=-15})).  

The set $\S(-25)$ comprises two equations of type 
I
and two equations of type 
II.
Three out of the four equations correspond to a zero rank elliptic curve $E_i$ at (\ref{e1e2}) for 
either $i=1$ or $i=2$.
Only the equation $5y_1^6+5y_2^6=2x^3$ remains. It is not difficult to prove that this 
equation is everywhere locally solvable. However, according to our general discussion 
at the beginning of section (\ref{generalities}),
a solution $(x,y_1,y_2)$ ($y_1y_2\neq 0$) gives rise to a point $Q_1=(10x/y_1^2,50y_1^3/y_2^3)$ 
on $E_1:=[0,0,0,0,-4\cdot 5^4]$. This elliptic curve has rank 1, with $E_1(\Q)$
generated by the point $P=(50,350)$ of infinite order. Let $Q_1=n\cdot P$. It is easily 
checked that 
\[ n\cdot P\equiv \O, (7,\pm 6), (42,\pm 6), (37,\pm 6) \pmod{43} \:.\]
Since the second coordinate of $Q_1$ equals $50y_1^3/y_2^3$, and the congruence 
$50Y^3\equiv\pm 6\!\!\pmod{43}$ is impossible, we conclude that the only possibility 
is $Q_1\equiv\O\!\!\pmod{43}$. This means that $y_2$ is divisible by 43, and consequently
$5y_1^6\equiv 2x^3\!\!\pmod{43}$, which is possible only if $y_1$ is divisible by 43; 
this contradicts the fact that $y_1,y_2$ are relatively prime.

The case $k=-29$ can be treated in complete analogy. The set $\S(-29)$ contains exactly 
three equations, two of which furnish a zero rank elliptic curve $E_i$. Only the equation 
$16y_1^6+29y_2^6=x^3$ remains, which is everywhere locally solvable. A solution 
$(x,y_1,y_2)$ gives a point $Q_2$ on $E_2:=[0,0,0,0,-16\cdot 29^3]$ of rank one and generator 
$P=(33085897/606^2, 129969272827/606^3)$. The relation
$Q_2=(29x/y_1^2,29y_2^3/y_1^3)=n\cdot P$ is proven impossible as above, working modulo 19.
%
\section{Application of Elliptic Chabauty}
\label{Chabauty} 
%
In the range $|k|\leq 50$ we are left with the following values
\begin{equation} \label{k-sieve 4}
k=-47,-39,-28,-15,-11,3,10,15,17,24,35,43,48\:.
\end{equation}
A natural approach to equation (\ref{xyeq}) is to factorize it over the field 
$K=\Q(\th)$, where $\th^3=k$. The maximal order $\O_K$ of $K$ has one fundamental unit 
$\ep(\th)$ which we normalize by sign to satisfy $\ep(k^{1/3})>0$ for the real
cube root $k^{1/3}$ of $k$. For $k \not \equiv \pm 1\!\! \pmod{9}$, the ideal 
$\langle 3 \rangle$ factors as $\langle3\rangle =\idp_3^3$; and when 
$k \equiv \pm 1 \!\!\pmod{9}$, then $\langle 3 \rangle =\idp_3 \idp_3'^2$.
We deduce from (\ref{xyeq}) the following ideal equations
\begin{equation} \label{quartquadideal}
\langle x^4-\th x^2 y^2+\th^2 y^4 \rangle =\idc\ida^2, \quad 
\langle x^2+\th y^2 \rangle  = \idc\idb^2, 
\quad \langle z \rangle =\idc\ida\idb
\end{equation}
for ideals $\ida$, $\idb$, $\idc$ of $\O_K$, where
\[ \idc= \langle x^4-\th x^2 y^2+\th^2 y^4\,,\, x^2+\th y^2 \rangle =
            \langle x^2+\th y^2 \,,\,3\th^2\rangle
\:.\]
If the highest power of every rational prime dividing $k$ is odd, then 
\[ \idc = \begin{cases}
          \O_K & \mbox{if $k\not\equiv 8 \!\!\pmod{9}$} \\
          \idp_3\idp_3' & \mbox{if $k\equiv 8 \!\!\pmod{9}$}
     \end{cases}
\]
If there exist rational primes whose highest power dividing $k$ is $2$ or $4$, then $\idc$ is as before, 
times an ideal $\idc_0$ which is divisible only by prime ideals over such rational primes
with exponents bounded by a small explicit integer.

In general, the equations at (\ref{quartquadideal}) are equivalent to
a system of element equations
\begin{equation}
\label{quartquadelement}
x^4-\th x^2 y^2+\th^2 y^4 = cw^2, \quad x^2+\th y^2 = c v^2, \quad z=cwv,
\end{equation}
for finitely many $c \in K$. We focus on the quartic curve $C$ defined by the first equation
at (\ref{quartquadelement}).
If this quartic represents an elliptic curve over $K$ (so in practice, if we can find 
a solution $(x,y,w)$ in $K$), then we seek points $(x/y,w/y^2)$ on the curve 
subject to the rationality condition $x/y \in \Q$.
We make extensive use of the Magma routines (inter alia)
{\tt PseudoMordellWeilGroup} to compute a subgroup of odd finite index in
$C(K)$, and {\tt Chabauty} for computing the $K$-points on $C$ with
prescribed rationality condition. Note that for {\tt Chabauty} a prime $p$ is
also needed as an input for the routine to work $p$-adically and, in certain
instances, an appropriate set of auxiliary primes. The choice of $p$ and auxiliary 
primes, if any, seems to be {\it ad exemplum}. 

As a characteristic example, consider $k=-28$.
Here, $K=\Q(\th)$ with $\th^3=-28$; the fundamental unit is $\ep=(-2+2\th+\th^2)/6$, 
the class number of $\O_K$ is 3 and we have four quartics $C_i$ corresponding to 
$c=1,\ep, 4+\th, \ep(4+\th)$, respectively.  Both quartics $C_2$ and $C_4$ are unsolvable 
at the prime ideal $\langle 2,(4+2\th+\th^2)/6\rangle$.  For $C_1$, Magma routines applied 
with {\tt Chabauty} using $p=11$ show that there are no solutions with $y\neq 0$. 
As for $C_3$, it has several points over $K$, the ``simplest'' one being
$(x/y,w/y^2)=(-2,6/(4+\th))$.
%
%
Using $p=5$ and auxiliary prime $17$,
{\tt Chabauty} reveals $x/y=\pm 2$ as the only possibility with $y\neq 0$, which returns 
$(\pm X, \pm Y)=(2,6)$ as the only finite points on (\ref{kcurve}). 
In an analogous manner, we solve equation (\ref{xyeq}) and hence (\ref{kcurve})
when $k=3$, $10$, $17$, $24$, $35$, and $48$;
see the table in section (\ref{Table}).
Note that these values of $k$ do not cover the full list (\ref{k-sieve 4}). Why not?
This is discussed in sections (\ref{difficultcases}) and (\ref{defeat}).
%
\section{The difficult cases $k=15$, $43$, $-11$, $-15$} \label{difficultcases}
%
%
We were unable to apply directly methods of the previous section for the values
$k=-47$, $-39$, $-15$, $-11$, $15$, $43$. In most cases, for at least one value of $c$ at
(\ref{quartquadelement}), our machines were unable to compute
the relevant Mordell-Weil groups over the cubic number field $K$. Oftentimes
the Selmer bound on the rank was equal to 3, but at most one non-torsion point could
be found, indicating the possible presence of a non-trivial Shafarevic-Tate group
(of course, if the rank is actually equal to 3, then elliptic Chabauty arguments
over a cubic number field must fail).  It is possible that a standard descent 
could be carried out by hand in such cases, but the calculations are rather daunting. 
Such an example of Selmer rank bound 3 occurs for instance when $k=15$ with 
$c=1-30\th+12\th^2$, $\th^3=15$ ($c$ being a fundamental unit in $\Q(\th)$).
In such cases we also investigated the quartic cover 
obtained by eliminating $x$ at (\ref{quartquadelement}), which results in the curve
\begin{equation}
\label{vyquartic}
c v^4- 3 \th v^2 y^2 +\frac{ 3\th^2}{c} y^4 = w^2.
\end{equation}
Here, and in several other cases, the $K$-rank of the curve (\ref{vyquartic}) could be
computed exactly but turned out to equal 3 (and in one case, 4), so not strictly less
than the degree of $\Q(\theta)$; hence {\tt Chabauty} is not applicable.\\ 
We adopt an alternative approach, which is successful for four of the values, 
namely $k=43$, $15$, $-11$, $-15$.  This approach involves factorization over an appropriate 
{\it quadratic} number field of the equation of type 
I, II, III, IV
that is causing the difficulty.
%
%
\subsection{Case $k=15$} \label{k=15}
%
%
Here, the finite points are $(\pm1,\pm4)$, $(\pm1/2,\pm31/8)$.\\ \\
We have to solve the equation $x^6+15y^6=z^2$, and following the approach of 
section (\ref{elementary}) we are left with solving
each of the two equations
\begin{equation}
\label{e1}
-5y_1^6 + 3y_2^6 = 2 x^3, \quad y_1,y_2 \mbox{ odd}, \quad (5y_1,3y_2)=1,
\end{equation}
and
\begin{equation}
\label{e2}
-16 y_1^6 +15 y_2^6 = x^3, \quad (2y_1,15y_2)=1.
\end{equation}
We work in $K=\Q(\phi)$, $\phi^2=15$, with maximal order $\O_K$. 
The ideal $\langle 2 \rangle =\idp_2^2$, $\idp_2=\langle 2,1+\phi \rangle$,
and the ideal classgroup of $\O_K$ is of order 2, 
generated by $\idp_2$. We have 
$\langle 3 \rangle=\idp_3^2$, 
$\idp_3=\langle 3,\phi \rangle$; $\langle 5 \rangle=\idp_5^2$, $\idp_5=\langle 5,\phi \rangle$. 
A fundamental unit is $\ep=4-\phi$.
\subsubsection{Equation (\ref{e1})}
We have
\[ (y_1^3 \phi+3 y_2^3) (y_1^3\phi-3 y_2^3) = -6 x^3, \]
and the gcd of the two factors on the left is precisely
$\idp_2 \idp_3 = \langle -3+\phi \rangle$.  Thus
\[ \langle y_1^3 \phi+3 y_2^3 \rangle = \langle -3+\phi \rangle \mathcal{A}^3, \]
for an ideal $\mathcal{A}$ prime to $\idp_2 \idp_3 \idp_5$ of $\O_K$.
Since $\mathcal{A}^3$ is principal and the class-number of $K$ is 2, 
then $\mathcal{A}$ is principal. So there exists
$y_3 \in \O_K$ satisfying
\begin{equation}
\label{e11}
y_1^3 \phi + 3 y_2^3 = \ep^i (-3+\phi) y_3^3, \quad i=0,\pm1.
\end{equation}
Equation (\ref{e11}) represents a curve of genus 1 over $K$,
which is locally unsolvable at $\idp_3$ when $i=1$; further, since
$\ep^{-1}(-3+\phi)=(3+\phi)$, the curve with $i=-1$ is simply the
conjugate of the curve with $i=0$. It suffices therefore to find
all points with rational $y_1:y_2$ on
\begin{equation}
\label{ee1}
y_1^3 \phi + 3 y_2^3 =(-3+\phi) y_3^3,
\end{equation}
which is an elliptic curve since it possesses the point
$(y_1,y_2,y_3)=(1,-1,1)$. The Magma routine {\tt PseudoMordellWeilGroup} tells us 
that the curve has rank 3 over $K$, however, exceeding the degree of
$K$, and so we cannot directly use the elliptic Chabauty method and Magma's
relevant routines.  We overcome this difficulty as follows.
At (\ref{ee1}), put $y_3=a+b\phi$, where $(a,15)=1$, $a+b \equiv 1 \bmod{2}$,
and where, changing the sign of each $y_i$ if necessary, we may assume $a \equiv 1 \bmod{3}$. Then
\begin{eqnarray} 
a^3 - 9a^2b + 45a b^2 - 45b^3 & = & y_1^3  \nonumber \\
-a^3 + 15a^2b - 45a b^2 + 75b^3 & = & y_2^3, \label{ec1}
\end{eqnarray}
defining rational elliptic curves of rank 1, 2, respectively.
The cubics at (\ref{ec1}) factor over the field $L=\Q(\psi)$, $\psi^3-3\psi+8=0$. 
A fundamental unit in the ring of integers $\O_L=\Z[1,\psi,\psi^2]$ is $\eta=5+2\psi$.
We have $\langle 2 \rangle =\idP_2 \idP_2'^2$ with $\idP_2=\langle 2,\psi \rangle$, 
$\idP_2'=\langle 2,1+\psi \rangle$;
$\langle 3 \rangle=\idP_3^3$ with $\idP_3 = \langle 3,2+\psi \rangle$; and
the classgroup is of order 3 generated by $\idP_2$, with $\idP_2^3=\langle \psi \rangle$.  \\
From (\ref{ec1})
\begin{eqnarray*}
(a + (\psi^2+\psi-5)b)(a^2 + (-\psi^2-\psi-4)ab + (3\psi^2-3\psi+9)b^2) & = & y_1^3, \\
(a + (-\psi^2+\psi-3)b)(a^2 + (\psi^2-\psi-12)ab + (-5\psi^2-5\psi+25)b^2) & = & -y_2^3.
\end{eqnarray*}
By the assumptions on $a,b$, it is straightforward to verify in each equation that the two 
factors on the left, considered as principal ideals, are coprime, and hence equal to ideal cubes.
Now an ideal equation $\langle u \rangle=\mathcal{B}^3$
implies one of the principal ideal equations
\[ \langle u \rangle= \langle v \rangle^3, \qquad \langle \psi^2 \rangle \langle u \rangle=\langle v \rangle^3, \qquad \langle \psi \rangle  \langle u \rangle= \langle v \rangle^3, \]
according as $\mathcal{B} \sim 1, \idP_2, \idP_2^2$.
So without loss of generality we deduce element equations
\begin{eqnarray}
\psi^{i_1}(a+(\psi^2+\psi-5)b) & = & \eta^{j_1} c_1^3, \nonumber \\
\psi^{i_2}(a^2+(-\psi^2-\psi-4)ab+(3\psi^2-3\psi+9)b^2) & = & \eta^{-j_1} c_2^3, \label{ab1}
\end{eqnarray}
with $i_1+i_2 \equiv 0 \bmod{3}$, $j_1=0,\pm 1$, and
\begin{eqnarray}
\psi^{i_3}(a+(-\psi^2+\psi-3)b) & = & \eta^{j_2} c_3^3, \nonumber \\
\psi^{i_4}(a^2+(\psi^2-\psi-12)ab+(-5\psi^2-5\psi+25)b^2) & = & \eta^{-j_2} c_4^3, \label{ab2}
\end{eqnarray}
with $i_3+i_4 \equiv 0 \bmod{3}$, $j_2=0,\pm 1$; and 
$c_i$, $i=1,..,4$, in $\O_L$. \\
We denote the $\idP_3$-adic (additive) valuation
of an element $\alpha \in \O_L$ by $\nu(\alpha)$, the highest power of $\idP_3$ 
dividing $\langle \alpha \rangle$, and for reference list here the valuations 
of certain elements of $\O_L$ (the coefficients of the polynomials 
occurring at equations (\ref{ab1}), (\ref{ab2})):
\[ \nu(\psi^2+\psi-5)=2; \qquad \nu(-\psi^2-\psi-4)=2; \qquad \nu(3\psi^2-3\psi+9)=4; \]
\[ \nu(-\psi^2+\psi-3)=1; \qquad \nu(\psi^2-\psi-12)=1; \qquad \nu(-5\psi^2-5\psi+25)=2; \]
further,
\[ \nu(\eta-1)=1; \qquad \nu(\psi-1)=1, \]
so that $\psi^3 \equiv 1 \bmod{3}$.
Thus (\ref{ab1}) and (\ref{ab2}) imply that $c_i^3 \equiv 1 \bmod{\idP_3}$, $i=1,..,4$,
so that $c_i \equiv 1 \bmod{\idP_3}$, whence $c_i^3 \equiv 1 \bmod{3}$, $i=1,..,4$.
We consider the first equations at (\ref{ab1}) and (\ref{ab2}) mod 3, distinguishing 
cases according as to the residue class of $b$ mod 3.\\
Subcase (i): $b \equiv 0 \bmod{3}$. Then
\[ \psi^{i_1} \equiv \eta^{j_1}, \qquad \psi^{i_3} \equiv \eta^{j_2},  \]
which forces $(i_1,i_3,j_1,j_2)=(0,0,0,0)$, and $(i_2,i_4)=(0,0)$.\\
Subcase (ii): $b \equiv 1 \bmod{3}$. Then
\[ \psi^{i_1}(\psi^2+\psi-4) = \eta^{j_1}, \qquad \psi^{i_3}(-\psi^2+\psi-2) = \eta^{j_2}, \]
which forces $(i_1,i_3,j_1,j_2)=(1,2,2,2)$ and $(i_2,i_4)=(2,1)$.\\
Subcase (iii): $b \equiv -1 \bmod{3}$. Then
\[ \psi^{i_1}(-\psi^2-\psi+6)) = \eta^{j_1}, \qquad \psi^{i_3}(\psi^2-\psi+4)) = \eta^{j_2}, \]
which forces $(i_1,i_3,j_1,j_2)=(2,0,1,2)$ and $(i_2,i_4)=(1,0)$.\\
Consequently, when we form the following equation using factors 
from (\ref{ab1}), (\ref{ab2}),
\[ \psi^{i_1+i_4}(a+(\psi^2+\psi-5)b)(a^2+(\psi^2-\psi-12)ab+(-5\psi^2-5\psi+25)b^2) = \eta^{j_1-j_2} c^3, \]
we have the three possibilities:
\begin{eqnarray*}
(a+(\psi^2+\psi-5)b)(a^2+(\psi^2-\psi-12)ab+(-5\psi^2-5\psi+25)b^2) & = & c^3, \\
\psi^2(a+(\psi^2+\psi-5)b)(a^2+(\psi^2-\psi-12)ab+(-5\psi^2-5\psi+25)b^2) & = & c^3, \\
\psi^2(a+(\psi^2+\psi-5)b)(a^2+(\psi^2-\psi-12)ab+(-5\psi^2-5\psi+25)b^2) & = & \eta^{-1} c^3,
\end{eqnarray*}
each representing an elliptic curve over $L$ having no $L$-torsion and of
rank 2 over $L$. Note that now the rank is less than the degree of L
which permits us to {\em attempt} the elliptic Chabauty method using
the Magma routines. And indeed, these routines (working $23$-adically
with auxiliary prime 53) show the following. The only point with 
rational $a:b$ on the first curve is when $a:b=1:0$, returning $(X,Y)=(1,4)$ 
on (\ref{kcurve}). 
Further, there are no such points on the second and
third curves, working $p$-adically respectively with $p=19$ and auxiliary prime
$13$, and $p=23$ with auxiliary primes $29,43$. \\
Remark: A gloss is needed for the Magma computations in this third case.
The routine {\tt Chabauty} assures us that there are no points on the third curve
with rational $a:b$, {\it provided} that the index of the group
of points (with generators $P_1$, $P_2$, say) output by 
{\tt PseudoMordellWeilGroup} in the full Mordell-Weil group is prime to 3.
This amounts to showing that none of $P_1$,$P_2$, $P_1 \pm P_2$ is divisible by 3
in the full Mordell-Weil group, and this is readily verified by Magma.
\subsubsection{Equation (\ref{e2})}
We have
\[ (y_2^3\phi +4 y_1^3) (y_2^3 \phi -4 y_1^3) = x^3, \]
and the two factors on the left are coprime, so that
\[ \langle y_2^3\phi +4 y_1^3 \rangle = \mathcal{A}^3, \]
for some ideal $\mathcal{A}$ prime to $\idp_2 \idp_3 \idp_5$ of $\O_K$.
Then as above, there exists $y_3 \in \O_K$ satisfying
\begin{equation}
\label{e21}
y_2^3\phi +4 y_1^3 = \ep^i y_3^3, \quad i=0, \pm1.
\end{equation}
Equation (\ref{e21}) is not locally solvable at $\idp_3$ when $i=0$, and        
the curves corresponding to $i=\pm1$ are conjugate; so                          
it suffices to consider only the case $i=1$:                                    
\begin{equation}                                                                
\label{ee2}                                                                     
y_2^3\phi +4 y_1^3 = (4-\phi) y_3^3,                                            
\end{equation}                                                                  
which is an elliptic curve since it possesses the point                         
$(y_1,y_2,y_3)=(1,-1,1)$. The $K$-rank of the curve is 3. \\ 
At (\ref{ee2}), put $y_3=a+b\phi$, with $(a,15)=1$, $a+b \equiv 1 \bmod{2}$, and where, without
loss of generality, $a \equiv 1 \bmod{3}$, to give
\begin{eqnarray}
-a^3 + 12a^2b - 45a b^2 + 60b^3 & = & y_2^3  \nonumber \\
4a^3 - 45a^2b + 180a b^2 - 225b^3 & = & 4 y_1^3, \label{ec2}
\end{eqnarray}
defining rational elliptic curves of rank 2, 1, respectively.
Note in this latter pair of equations that the second equation implies
$b(a^2+b^2) \equiv 0 \bmod{4}$, so that necessarily $a$ is odd,
$b \equiv 0 \bmod{4}$. \\
From (\ref{ec2}), factoring over $L$:
\begin{eqnarray*} 
(a + (\psi-4)b)(a^2 + (-\psi-8)ab + (\psi^2+4\psi+13)b^2) & = & -y_2^3 \\
(a + (-\psi^2-4\psi-13)b/4)(a^2 + (\psi^2+4\psi-32)ab/4 + 15(-\psi+4)b^2/4) & = & y_1^3
\end{eqnarray*} 
and again in each equation the two factors on the left, considered as principal
ideals, are coprime and hence ideal cubes.
Just as above, we deduce element equations 
\begin{eqnarray}
\psi^{i_1}(a+(\psi-4)b) & = & \eta^{j_1} c_1^3, \nonumber \\
\psi^{i_2}(a^2 + (-\psi-8)ab + (\psi^2+4\psi+13)b^2) & = & \eta^{-j_1} c_2^3, \label{ab5}
\end{eqnarray}
with $i_1+i_2 \equiv 0 \bmod{3}$, $j_1=0,\pm 1$, and
\begin{eqnarray}
\psi^{i_3}(a + (-\psi^2-4\psi-13)b/4) & = & \eta^{j_2} c_3^3, \nonumber \\
\psi^{i_4}(a^2 + (\psi^2+4\psi-32)ab/4 + 15(-\psi+4)b^2/4) & = & \eta^{-j_2} c_4^3, \label{ab6}
\end{eqnarray}
with $i_3+i_4 \equiv 0 \bmod{3}$, $j_2=0,\pm 1$ ; and $c_i$, $i=1,..,4$, 
in $\O_L$.\\
We have the following valuations:
\[ \nu(\psi-4)=1; \qquad \nu(-\psi-8)=1; \qquad \nu(\psi^2+4\psi+13)=2; \]
\[ \nu(\psi^2+4\psi-32)=2; \qquad \nu(15(-\psi+4))=4, \]
so that as before, $c_i^3 \equiv 1 \bmod{3}$. We consider the first equations 
at (\ref{ab5}) and (\ref{ab6}) mod 3, distinguishing cases according as to the 
residue class of $b$ mod 3.\\
Subcase (i): $b \equiv 0 \bmod{3}$. Then
\[ \psi^{i_1} \equiv \eta^{j_1}, \qquad \psi^{i_3} \equiv \eta^{j_2}, \] 
forcing $(i_1,i_3,j_1,j_2)=(0,0,0,0)$ and $(i_2,i_4)=(0,0)$.\\
Subcase (ii): $b \equiv 1 \bmod{3}$. Then
\[ \psi^{i_1}(\psi-3) = \eta^{j_1}, \qquad \psi^{i_3}(-\psi^2-4\psi-9) = \eta^{j_2}, \]
forcing $(i_1,i_3,j_1,j_2)=(2,2,0,1)$ and $(i_2,i_4)=(1,1)$.\\
Subcase (iii): $b \equiv -1 \bmod{3}$. Then
\[ \psi^{i_1}(-\psi+5) = \eta^{j_1}, \qquad \psi^{i_3}(\psi^2+4\psi+17)) = \eta^{j_2}, \]
forcing $(i_1,i_3,j_1,j_2)=(0,1,1,2)$ and $(i_2,i_4)=(0,2)$.\\
Consequently, when we form the following equation using factors 
from (\ref{ab5}), (\ref{ab6}),
\[ \psi^{i_1+i_4}(a+(\psi-4)b)(a^2 + (\psi^2+4\psi-32)ab/4 + 15(-\psi+4)b^2/4) = \eta^{j_1-j_2} c^3, \]
we have the three possibilities:
\begin{eqnarray*}
(a+(\psi-4)b)(a^2 + (\psi^2+4\psi-32)ab/4 + 15(-\psi+4)b^2/4) & = & c^3, \\
(a+(\psi-4)b)(a^2 + (\psi^2+4\psi-32)ab/4 + 15(-\psi+4)b^2/4) & = & \eta^{-1} c^3, \\
\psi^2(a+(\psi-4)b)(a^2 + (\psi^2+4\psi-32)ab/4 + 15(-\psi+4)b^2/4) & = & \eta^{-1} c^3,
\end{eqnarray*}
again each representing an elliptic curve of rank $2$ over $L$.
Working $13$-adically, having shown that the group index is prime to $3$ (actually
we need it prime to $6$, but the construction of the routine 
{\tt PseudoMordellWeilGroup} guarantees that the index is odd), 
we find that the first curve has a point with rational $a:b$ only
at $a:b=1:0$, returning $(X,Y)=(-\frac{1}{2},\frac{31}{8})$.
Working $13$-adically, with auxiliary prime $43$, we find there are no points 
with rational $a:b$ on the second curve (where we needed to show the group 
index is prime to $3$);
and working $13$-adically, with auxiliary prime 17, we find there
are no points with rational $a:b$ on the third curve (where we needed to
show the group index is prime to 3).
In summary, the only finite points 
on (\ref{kcurve}) are given by $(\pm X, \pm Y)=(1,4)$, $(\frac{1}{2}, \frac{31}{8})$.
%
%
\subsection{Case $k=43$} \label{k=43}
%
%
Here, the finite points are $(\pm3/2,\pm59/9)$, $(\pm7/3,\pm386/27)$. \\ \\
This case is similar to the previous one.  The set $S(43)$ derived in section
(\ref{elementary args}) contains the three globally solvable equations:
\begin{equation}
\label{S(43)}
2x^3=-y_1^6+43 y_2^6, \qquad x^3=-16 y_1^6+43y_2^6, \qquad x^3=-688 y_1^6+y_2^6.
\end{equation}
We factor the first two equations over the field $\Q(\phi)$, where $\phi^2=43$.
As before, we are led to working in a cubic field $\Q(\psi)$,
where now $\psi^3-\psi^2+\psi-9=0$. In analogy to the equations at the end of 
sections (4.1.1) and (4.1.2), we obtain the equations:
\begin{align*}
((-\psi^2 - \psi& - 3)a + (6\psi^2 + 6\psi + 17)b) \\
((-10\psi^2 - 13\psi& - 34)a^2 + (124\psi^2 + 163\psi + 417)a b + (-387\psi^2 - 516\psi - 1290)b^2) \\
& = \ep^l (-\psi^2+\psi+3) (2-\psi)^2 c^3,
\end{align*}
where $l=0,1,2$, and
\begin{align}
\label{second43}
((5\psi^2+6\psi+20)a& + \frac{1}{2}(-\psi^2-2\psi-5)b) \\ \nonumber
(\frac{1}{2}(-\psi^2-4\psi&-5)a^2 + \frac{1}{2}(\psi^2-2\psi+1)ab - b^2) = \frac{1}{2}(-5\psi^2 + 4\psi + 17) \ep^l c^3, 
\end{align}
where $l=0,1,2$.
Just as before, these equations are amenable to the computer routines, and deliver precisely the two known finite solutions.
It is worth mentioning that the relevant Mordell-Weil group
on this latter curve of rank 2 over $\Q(\psi)$ could not be computed directly with Magma because of 
the large height of one of the generators: the generators on (\ref{second43}) may be taken as $(a,b,c)=(0,-7\psi+16,1)$ and
\begin{align*}
(\frac{1}{2}(-24963589\psi^2 + 48373018\psi + 20008291), \quad 18015880\psi^2 - 132297936\psi + 209249791,& \\
31418694\psi^2 + 46376736\psi + 127224108).&
\end{align*}
It was necessary to consider a 3-isogenous curve where
by luck the relevant group could be computed directly. This gives rise to a full rank
subgroup on (\ref{second43}), successfully feeding into the {\tt Chabauty} routine. \\
For the third equation at (\ref{S(43)}), the solution is accomplished by the elementary methods of
section (\ref{elementary args}).
%
%
\subsection{Case $k=-11$} \label{k=-11}
%
%
Here, the finite points are $(\pm3/2,\pm5/8)$. \\ \\ 
We have to solve $x^6-11y^6=z^2$, where $x,11y,z$ are pairwise relatively prime.
Following the approach of (\ref{elementary}), we are left with solving
\begin{equation} \label{minus11.11}
 x^3 = 16 y_1^6 + 11 y_2^6,
\end{equation}
where $y=2y_1y_2$ and $(2y_1,11y_2)=1$. Considering equation (\ref{minus11.11}) $\bmod{\,9}$
shows that, if either $y_1$ or $y_2$ is divisible by 3, then both $y_1,y_2$ are divisible by 3, a contradiction;
hence $y_1y_2\not\equiv 0\pmod{3}$. We factorize equation (\ref{minus11.11}) as
$(4y_1^3+\sqrt{-11}y_2^3)(4y_1^3-\sqrt{-11}y_2^3)=x^3$, where the two factors
on the left hand side are coprime in the ring of integers of $\Q(\sqrt{-11})$,
a field of class number 1.  It follows that 
$4y_1^3+\sqrt{-11}y_2^3=(a+b\frac{(1+\sqrt{-11})}{2})^3$ with $a,b\in\Z$,
and hence that
\begin{equation} \label{k11.initial}
b(3a^2-3ab-2b^2)= 2y_2^3, \qquad (2a-b)(a^2-ab-8b^2)=(2y_1)^3
\end{equation}
with $(a,b)=1$. These equations define elliptic curves of ranks 2, 1, respectively.  
Note that the second equation $\bmod{\,9}$ shows that $b\not\equiv 0\bmod{3}$. \\
Since $(b,3a)=1$, the two factors on the left hand side of the first equation at
(\ref{k11.initial}) are coprime.
Also, $(2a-b,a^2-ab-8b^2)=(2a-b, 33)=1$, since $(y_1,33)=1$. It follows that
\[ b=2^i\cdot\mbox{cube}, \quad 3 a^2 - 3 a b - 2 b^2=2^{1-i}\cdot\mbox{cube}, \;\; i=0,1, \]
\[ 2a-b = \mbox{cube}, \quad a^2-a b-8 b^2 = \mbox{cube}. \]
If $i=1$, then we deduce
\[ b=2\be^3, \quad 2a-b=8\al^3, \quad a^2-a b-8 b^2 = \ga^3 , \]
so that
\[ 16 \al^6 - 33 \be^6 = \ga^3. \]
This is an equation of type (\ref{ABsextic}), and taking $(A,B)=(16,-33)$, we discover that
the corresponding elliptic curve $E_3$ at (\ref{e3}) has rank 0; and no solutions arise
for $a$, $b$. Thus $i=0$, and we have
\begin{equation} \label{cubeeqs}
b =\mbox{cube}, \quad 2a-b = \mbox{cube}, \quad a^2-a b-8 b^2 = \mbox{cube}, \quad 3 a^2 - 3 a b - 2 b^2=2\cdot\mbox{cube}.
\end{equation}
Note that, in the above equations, $a^2-ab-8b^2$ is not divisible by 11 because it is a factor
of $y_1$. Also, $3a^2-3ab-2b^2$ is not divisible by 11 for, otherwise, it would be divisible
by $11^2$ which implies $b\equiv 0\pmod{11}$, hence also $a\equiv 0\pmod{11}$; a contradiction.
These observations will be needed below, when we calculate greatest common divisors.

We work in the field $\Q(\xi)$, $\xi^2-\xi-8=0$. The class-number is 1, a fundamental
unit is $\ep=19+8\xi$, and we have the prime factorization $2=(2+\xi)(-3+\xi)$.

The latter two equations at (\ref{cubeeqs}) may be written as follows:
\[ (a-\xi b)\;(a+(-1+\xi)b) = \mbox{cube}, \quad 
((5+2\xi)a+(2+\xi)b)\;((7-2\xi)a+(3-\xi)b) = 2 \cdot \mbox{cube}. \]
The greatest common divisor $(a-\xi b, a+(-1+\xi)b)=(a-\xi b, 1-2\xi)=(a-\xi b, \sqrt{33})=1$,
since  $(y_2,33)=1$. Thus
\[ a-\xi b = \ep^j\;\mbox{cube}, \quad a+(-1+\xi)b = \ep^{-j}\;\mbox{cube}, \quad j=0,\pm1.  \]
Further, the greatest common divisor 
$((5+2\xi)a+(2+\xi)b, (7-2\xi)a+(3-\xi)b)= ((5+2\xi)a+(2+\xi)b, 1-2\xi)=1$, as before.
Thus 
\[ (5+2\xi)a+(2+\xi)b = \ep^k\;\pi\;\mbox{cube}, \quad (7-2\xi)a+(3-\xi)b = \ep^{-k}\;\bar{\pi}\;\mbox{cube}, \quad k=0,\pm1, \]
for $\pi$, $\bar{\pi}$ equal to the two prime factors of 2.

Summing up, we have
\[ b = \mbox{cube}, \quad 2a-b=\mbox{cube}, \quad a+(-1+\xi)b = \ep^{-j}\;\mbox{cube}, \quad (5+2\xi)a+(2+\xi)b = \ep^k\;\pi\;\mbox{cube}, \]
for $\pi=2+\xi$ or $-3+\xi$.
We form the elliptic cubic
\[ b\;(a + (-1+\xi)b)\;((5+2\xi)a+(2+\xi)b) = \ep^l\;\pi\;\mbox{cube}, \quad l=-j+k, \]
where without loss of generality, $l=0,\pm1$.
Consider first $\pi=2+\xi$. For each value of $l=0,\pm1$, the corresponding elliptic curve
has rank 1 over $\Q(\xi)$, and we can apply the elliptic Chabauty method.
For $l=0$, with {\tt Chabauty} working 7-adically, the only solutions are $a:b=1:0\,,-7:3$ 
which are rejected in view of (\ref{k11.initial}).  For $l=1$, working $7$-adically, the 
only solutions are $a:b=1:0,\,0:1$ and only the second satisfies (\ref{k11.initial}), 
giving $y_1=1,y_2=-1$, hence $y=-2,x=3,z=\pm 5$, which returns the points
$(\pm X, \pm Y)=(\frac{3}{2}, \frac{5}{8})$ on (\ref{kcurve}).
For $l=-1$, working $19$-adically, the only solution is $a:b=1:0$ which is rejected
in view of (\ref{k11.initial}). \\
Second, take $\pi=3-\xi$. Again, each $l=0,\pm1$ results in a rank 1 elliptic curve.
If $l=0$, working 7-adically shows the only solutions to be $a:b=1:0,\,1:1$.
In view of (\ref{k11.initial}) the first is rejected, and the second gives
$y_1=-1,y_2=-1$, hence $y=2,x=3,z=\pm 5$, again returning the points
$(\pm X, \pm Y)=(\frac{3}{2},\frac{5}{8})$ on (\ref{kcurve}).
If $l=1$, working $7$-adically, the only solutions are $a:b=1:0,\,10:3$ which by
(\ref{k11.initial}) we reject.
If $l=-1$, working $19$-adically, the only solutions are $a:b=1:0$, which
by (\ref{k11.initial}) we reject. 
%
%
\subsection{Case $k=-15$} \label{k=-15}
%
%
Here, the finite points are $(\pm2,\pm7)$.\\ \\
We have to solve $x^6-15y^6=z^2$, where $x,15y,z$ are pairwise relatively prime.
Following the approach of section (\ref{elementary}), we are left with solving the
pair of equations
\begin{equation} \label{minus15.3}
80 y_1^6+3 y_2^6 = x^3\,,\quad \quad (10y_1,3y_2)=1\,,y=2y_1y_2, 
\end{equation}
and
\begin{equation} \label{minus15.5}
   y_1^6+15y_2^6 = 2x^3\,,\quad (y_1,15y_2)=1\,,\quad \mbox{$y_1y_2$ odd}\:.
\end{equation}
To deal with (\ref{minus15.3}), we use ideas of section (\ref{congruences}).
On the corresponding curve
\begin{equation} \label{minus15.1}
E_1: Y^2=X^3-375, \qquad (X,Y)=(\frac{5x}{y_2^2},\;\frac{100y_1^3}{y_2^3}) \:,
\end{equation}
the torsion is trivial, the rank is 1, and a generator $P$ is given
by $P=(10,25)$. 
We check that if $n\equiv 1,2,4,5 \bmod{6}$ then the $Y$-coordinate of $n\cdot P$
has odd numerator and denominator. Indeed, for $n=1,2,4,5$ this is straightforward; 
further, a symbolic computation shows that if we add to $6\cdot P$ a point 
$(\frac{u}{t^2}, \frac{v}{t^3})$, where $u,v,t$ are integers with $v t$ odd,
then the resulting point has $Y$-coordinate with odd numerator and denominator.
Hence for $n \equiv 1,2,5,6 \bmod{6}$, the $Y$-coordinate of $n\cdot P$ cannot 
have the required shape $100 y_1^3/y_2^3$ with $y_2$ odd. 
It remains to check the cases $n \equiv 0,3 \bmod{6}$.
But if $n$ is a multiple of 3, the
following are the possibilities for the $Y$-coordinate of $n\cdot P$: either
it is congruent to $\pm 9 \bmod{19}$, or its denominator is divisible by 19.
The first alternative is impossible because it
implies that $100(y_1/y_2)^3\equiv \pm 9 \bmod{19}$; and the second implies
$y_2\equiv 0 \bmod{19}$, but then the initial equation $80y_1^6+3y_2^6=x^3$ 
is impossible $\bmod{\,19}$ when $(y_1,y_2)=1$.\\
Now focus attention on equation (\ref{minus15.5}).
We work in the field $\Q(\th)$, where
$\th^2+\th+4=0$, with class-number $2$ and integral basis $1,\th$.
We have the ideal factorizations
\[ \langle 2\rangle =\idp_2\idp_2' \,,
      \quad \idp_2=\langle 2,1+\th\rangle\,, \quad \idp_2'=\langle 2,2+\th\rangle  \:.\]
The ideal-class $\idp_2$ generates the classgroup, and
$\idp_2^2=\langle 1+\th\rangle $.
Choosing appropriate signs for $y_1,y_2$, then
the ideal factorization of (\ref{minus15.5}) implies without loss of generality
\[ \langle y_1^3+(2\th+1)y_2^3\rangle =\idp_2\ida^3 \]
for some integral ideal $\ida$ such that $\idp_2\ida$ is principal.
This results in an equation
\begin{equation}
\label{minus15.magma}
 (\th+1)y_1^3+(\th-7)y_2^3=y_3^3
\end{equation}
for some $y_3 \in \Z[\th]$, which represents an elliptic curve over $\Q(\th)$ 
(note that it contains the point $(1,-1,2)$).
The $\Q(\th)$ rank is $1$, and Magma routines working
$17$-adically show that $(y_1,y_2,y_3)=(1,-1,2)$ is the only point
over $\Q(\th)$ with rational $y_1:y_2$. This gives $y=-1$,
and by (\ref{minus15.5}), $x=2$. Returning to (\ref{kcurve}), the only
finite points are $(\pm X, \pm Y)=(2,7)$.
\\
Remark: again, we needed to verify that the group of points output
by {\tt PseudoMordellWeilGroup} has index in the full Mordell-Weil group which is prime to 6.
%
%
\section{The unsolved equations} \label{defeat}
%
%
In the considered range of $k$, we are left with $k=-47,-39$. If we try
to apply the ideas of section (\ref{Chabauty}), then relevant Mordell-Weil groups
could not be computed. For example, at $k=-47$, $c=1$, the quartic curve 
at (\ref{quartquadelement}) has Selmer rank 3, with only one point
of infinite order found. The curve at (\ref{vyquartic}) has $K$-rank $4$.
Trying to apply the ideas of section (\ref{difficultcases}) for $k=-39$
leads to curves with bound on the rank 2, but where we are unable to find any points;
and similar obstructions arise for 
$k=-47$. We have tried various further attacks on these equations,
so far without success. It has been suggested to us that the computations
of this paper be automated to extend the calculations for $k$ in a range 
``somewhere in the thousands"; but without a mechanized $2$-descent
algorithm for elliptic curves over number fields, at the very least, even 
a range into the hundreds is well beyond our abilities. In exploring the unsolved 
cases of this section, we have resorted to much manual intervention in Magma
programming, primarily choosing appropriate models for curves and their isogenies
to replace the ones returned by the routine {\tt EllipticCurve}, which can have 
huge coefficients, possibly greatly increasing the running
time of the algorithms. When the number field has class-number exceeding 1, we know
of no uniform method for choosing models that are potentially better suited as input
to {\tt PseudoMordellWeilGroup}.\\ \\
For interest we searched the curves (\ref{kcurve}) for rational
points with height at most $20000$ in the range $|k| < 250000$. The maximum number
of points found was 22, at $k=1025$ with the finite points $(\pm X,\pm Y)=(2,33)$, 
$(1/4,2049/64)$, $(5/2,285/8)$, $(8,513)$, $(20/91,24126045/91^3)$; and
at $k=110160$, with finite points (all integral) $(\pm X,\pm Y)=(2,332)$, 
$(3,333)$, $(6,396)$, $(9,801)$, $(14,2764)$. The point of largest height that we 
observed in this range occurs for $k=-212860$, with point $(3866/427, 45259682826/427^3)$.
Increasing the search range and decreasing the height bound finds the curve
at $k=7547408$ with 26 points, the finite points occurring at $(\pm X, \pm Y)=$
\[ (4,2748), \; (\frac{1}{6}, \; \frac{593407}{6^3}), \; (\frac{7}{5},\frac{343407}{5^3}), \; (16,4932), \; (\frac{28}{3},\frac{77356}{3^3}),  \; (\frac{139}{10},\frac{3841869}{10^3}). \]
Note that the curve $y^2=x^6+(\frac{1}{4}a^{12}+1)$
automatically contains the points 
$(\pm X,\pm Y)=(\frac{1}{a^2},\frac{1}{a^2}+\frac{1}{2}a^6)$, 
$(a,1+\frac{1}{2}a^6)$, and $(\frac{1}{2}a^4, \frac{1}{8}a^{12}+1)$.
For the curve at $k=2089$ only 14 points were found,
but the (finite) points comprise $(\pm X,\pm Y)=(96/11,886825/11^3)$,
$(162/85,28389097/85^3)$, and $(289/90,41143681/90^3)$, remarkable for
their large height.

%
%
%
%
%
%
\section{All rational solutions to (\ref{kcurve}) in the range $|k|\leq 50$} \label{Table}
\begin{center}
All rational solutions $(\pm X,\pm Y)$ to $Y^2=X^6+k$,\;$1 \leq |k| \leq 50$
\\[1mm]
\begin{tabular}{|r|c|c|c|r|c|c|} \hline\hline
$k$ & Solutions $\pm(X,Y)$ & Ref.~section & & $k$ & Solutions $\pm(X,Y)$ & Ref.~section 
         \\ \hline\hline 
-50 & $\emptyset$ & (\ref{zero rank cE}) &  & 
 50 & $\emptyset$ & (\ref{zero rank cE}) \\
-49 & $\emptyset$ & (\ref{elementary args}) &  & 
 49 & (0,7) & (\ref{zero rank cE})  \\
-48 & (2,4) & (\ref{elementary args}) &  & 
 48 & (1,7) & (\ref{Chabauty})  \\ 
-47 & $(\frac{63}{10},\frac{249953}{1000})$ & (unproven) (\ref{defeat}) & & 
 47 & $\emptyset$ & (\ref{elementary args}) \\ 
-46 & $\emptyset$ & (\ref{zero rank cE}) & & 
46 & $\emptyset$ & (\ref{elementary args}) \\
-45 & $\emptyset$ & (\ref{elementary args}) & & 
45 & $\emptyset$ & (\ref{zero rank cE}) \\
-44 & $\emptyset$ & (\ref{zero rank cE}) & & 
44 & $\emptyset$ & (\ref{zero rank cE}) \\
-43 & $\emptyset$ & (\ref{elementary}) & & 
43 & $(\frac{3}{2},\frac{59}{9}),(\frac{7}{3},\frac{386}{27})$ & (\ref{k=43}) \\
-42 & $\emptyset$ & (\ref{zero rank cE}) & & 
42 & $\emptyset$ & (\ref{zero rank cE}) \\
-41 & $\emptyset$ & (\ref{zero rank cE}) & & 
41 & $\emptyset$ & (\ref{zero rank cE}) \\
-40 & $\emptyset$ & (\ref{zero rank cE}) & & 
40 & $\emptyset$ & (\ref{zero rank cE}) \\
-39 & $(2,5)$ & (unproven) (\ref{defeat}) & & 
39 & $\emptyset$ & (\ref{elementary args}) \\
-38 & $\emptyset$ & (\ref{zero rank cE}) & & 
38 & $\emptyset$ & (\ref{zero rank cE}) \\
-37 & $\emptyset$ & (\ref{zero rank cE}) & & 
37 & $\emptyset$ & (\ref{zero rank cE}) \\
-36 & $\emptyset$ & (\ref{zero rank cE}) & & 
36 & $(0,6),(2,10)$ & (\ref{elementary args}) \\
-35 & $\emptyset$ & (\ref{elementary}) & & 
35 & $(1,6)$ & (\ref{Chabauty}) \\
-34 & $\emptyset$ & (\ref{zero rank cE}) & & 
34 & $\emptyset$ & (\ref{zero rank cE}) \\
-33 & $\emptyset$ & (\ref{zero rank cE}) & & 
33 & $\emptyset$ & (\ref{zero rank cE}) \\
-32 & $\emptyset$ & (\ref{zero rank cE}) & & 
32 & $\emptyset$ & (\ref{zero rank cE}) \\
-31 & $\emptyset$ & (\ref{zero rank cE}) & & 
31 & $\emptyset$ & (\ref{elementary}) \\
-30 & $\emptyset$ & (\ref{zero rank cE}) & & 
30 & $\emptyset$ & (\ref{zero rank cE}) \\
-29 & $\emptyset$ & (\ref{congruences}) & & 
29 & $\emptyset$ & (\ref{zero rank cE}) \\
-28 & $(2,6)$ & (\ref{Chabauty}) & & 
28 & $\emptyset$ & (\ref{elementary args}) \\
-27 & $\emptyset$ & (\ref{zero rank cE}) & & 
27 & $\emptyset$ & (\ref{zero rank cE}) \\
-26 & $\emptyset$ & (\ref{zero rank cE}) & & 
26 & $\emptyset$ & (\ref{zero rank cE}) \\
-25 & $\emptyset$ & (\ref{congruences}) & & 
25 & $(0,5)$ & (\ref{zero rank cE}) \\
-24 & $\emptyset$ & (\ref{zero rank cE}) & & 
24 & $(1,5),(\frac{5}{2},\frac{131}{8})$ & (\ref{Chabauty}) \\
-23 & $\emptyset$ & (\ref{zero rank cE}) & & 
23 & $\emptyset$ & (\ref{zero rank cE}) \\
-22 & $\emptyset$ & (\ref{zero rank cE}) & & 
22 & $\emptyset$ & (\ref{zero rank cE}) \\
-21 & $\emptyset$ & (\ref{elementary}) & & 
21 & $\emptyset$ & (\ref{zero rank cE}) \\
-20 & $\emptyset$ & (\ref{zero rank cE}) & & 
20 & $\emptyset$ & (\ref{zero rank cE}) \\
-19 & $\emptyset$ & (\ref{zero rank cE}) & & 
19 & $\emptyset$ & (\ref{zero rank cE}) \\
-18 & $\emptyset$ & (\ref{zero rank cE}) & & 
18 & $\emptyset$ & (\ref{zero rank cE}) \\
-17 & $\emptyset$ & (\ref{zero rank cE}) & & 
17 & $(2,9),(\frac{1}{2},\frac{33}{8})$ & (\ref{Chabauty}) \\
-16 & $\emptyset$ & (\ref{zero rank cE}) & & 
16 & $(0,4)$ & (\ref{zero rank cE}) \\
\cline{1-3}\cline{5-7}
\end{tabular}
\end{center}

\pagebreak
\begin{center}
All rational solutions $(\pm X,\pm Y)$ to $Y^2=X^6+k$,\;$1 \leq |k| \leq 50$
\\ \footnotesize{(continued from previous page)}
\\[1mm]
\begin{tabular}{|r|c|c|c|r|c|c|} \hline\hline
$k$ & Solutions $\pm(X,Y)$ & Ref.~section & & $k$ & Solutions $\pm(X,Y)$ & Ref.~section 
         \\ \hline\hline 
-15 & $(2,7)$ & (\ref{k=-15}) & & 
15 & $(1,4)$ & (\ref{k=15}) \\
-14 & $\emptyset$ & (\ref{zero rank cE}) & & 
14 & $\emptyset$ & (\ref{zero rank cE}) \\
-13 & $\emptyset$ & (\ref{elementary args}) & & 
13 & $\emptyset$ & (\ref{zero rank cE}) \\
-12 & $\emptyset$ & (\ref{zero rank cE}) & & 
12 & $\emptyset$ & (\ref{zero rank cE}) \\
-11 & $(\frac{3}{2},\frac{5}{8})$ & (\ref{k=-11}) & & 
11 & $\emptyset$ & (\ref{elementary args}) \\
-10 & $\emptyset$ & (\ref{zero rank cE}) & & 
10 & $(\frac{3}{2},\frac{37}{8})$ & (\ref{Chabauty}) \\
-9 & $\emptyset$ & (\ref{zero rank cE}) & & 
9 & $(0,3)$ & (\ref{zero rank cE}) \\
-8 & $\emptyset$ & (\ref{zero rank cE}) & & 
8 & $(1,3)$ & (\ref{zero rank cE}) \\
-7 & $\emptyset$ & (\ref{zero rank cE}) & & 
7 & $\emptyset$ & (\ref{zero rank cE}) \\
-6 & $\emptyset$ & (\ref{zero rank cE}) & & 
6 & $\emptyset$ & (\ref{zero rank cE}) \\
-5 & $\emptyset$ & (\ref{zero rank cE}) & & 
5 & $\emptyset$ & (\ref{zero rank cE}) \\
-4 & $\emptyset$ & (\ref{zero rank cE}) & & 
4 & $(0,2)$ & (\ref{zero rank cE}) \\
-3 & $\emptyset$ & (\ref{zero rank cE}) & & 
3 & $\emptyset$ & (\ref{Chabauty}) \\
-2 & $\emptyset$ & (\ref{zero rank cE}) & & 
2 & $\emptyset$ & (\ref{zero rank cE}) \\
-1 & $(1,0)$ & (\ref{zero rank cE}) & & 
1 & $(0,1)$ & (\ref{zero rank cE}) \\
\cline{1-3}\cline{5-7}
\end{tabular}
\end{center}

\end{document}